\documentclass[12pt]{amsart}
\usepackage{amssymb}
\usepackage{amsfonts}
\usepackage{fullpage}
\usepackage{eufrak}
\input xypic
\xyoption{all}

\newtheorem{prop}{Proposition}
\newtheorem{lem}{Lemma}
\newtheorem{thm}{Theorem}
\newtheorem{cor}{Corollary}

\def\RR{{\mathbb R}}
\def\CC{{\mathbb C}}

\def\QQ{{\mathbb Q}}
\def\ZZ{{\mathbb Z}}

\def\cM{{\mathcal M}}

\def\cO{{\mathcal O}}

\def\cH{{\mathcal H}}
\title{Veech surfaces and complete periodicity in genus $2$} 
\author{Kariane Calta}
\begin{document}
\maketitle
\begin{abstract}
We announce several results pertaining to Veech surfaces and completely periodic translation surfaces in genus $2$.  A translation surface is a pair $(M, \omega)$ where $M$ is a Riemann
surface and $\omega $ is an abelian differential on $M$.  Equivalently, a
translation surface is a two-manifold which has transition functions which
are translations and a finite number of conical singularities arising from
the zeroes of $\omega$.  A saddle connection is a geodesic segment connecting two (not necessarily distinct) conical singularities with no singularities in its interior.

A direction $v$ on a translation surface is completely periodic if
any trajectory in the direction $v$ is either closed or ends in a
singularity, i.e. if the surface decomposes as a union of cylinders in the
direction $v$.  Then, we say that a translation surface is completely periodic if any direction in which there is a saddle connection joining a singularity to itself or at least one cylinder of closed trajectories
is completely periodic. There is an action of the group $SL(2, \RR)$ on the
space of translation surfaces.  A surface which has a lattice stabilizer
under this action is said to be Veech.  Veech proved that any Veech
surface is completely periodic, but the converse is false. 

In this announcement, we use the "$J$-invariant" of Kenyon and Smillie to obtain a classification of all Veech surfaces in the space $\cH(2)$ of genus $2$
translation surfaces with corresponding abelian differentials which have a
single double zero.  Furthermore, we obtain a classification of all
completely periodic surfaces in genus $2$. 
 \end{abstract}

A translation surface is a pair $(M, \omega)$ where $M$ is a closed Riemann surface and $\omega$ is an abelian differential on $M$.  Away from the zeroes of $\omega$, a chart $z$ can be chosen so that $\omega = dz$, which determines a Euclidean metric in that chart.  The change of coordinates away from the zeroes of $\omega$ are of the form $z \to z + c$.  In the neighborhood of a zero, one can choose a coordinate $z$ so that $\omega = z^k dz$.  The total angle about this zero is $2(k+1)\pi$ and the order of the zero is $k$.  Such zero is called a conical singularity. As a result of this definition, a translation surface can also be thought of as a two-manifold with transition functions which are translations and which has a finite number of conical singularities, arising from the zeroes of $\omega$. Each of these singularities has total angle $2\pi n$ for $n \in \ZZ^+$.  Thus, one can also think of a translation surface as a finite number of polygons in the plane $\RR^2$, glued along parallel sides.  Throughout, we will often denote a translation surface by $S$, and this notation should be thought of as incorporating both the Riemann surface $M$ and the differential $\omega$. 

If $S$ is a translation surface, a {\em saddle connection} is a geodesic connecting two (not necessarily distinct) conical singularities of $S$ which has no singularities in its interior.  Note that since the metric on $S$ is Euclidean away from the zeroes of $\omega$, saddle connections are straight lines, but in general, geodesics are unions of saddle connections where the line segments comprising the geodesic are allowed to switch directions at the conical singularities of $S$.  A closed geodesic on $S$ determines  a cylinder of parallel freely homotopic closed curves such that the boundaries of the cylinder consist of unions of saddle connections. A direction $v$ on $S$ is said to be {\em completely periodic} if all the trajectories in
the direction $v$ are either closed or terminate at a singularity. In
other words, $v$ is completely periodic if in the direction $v$, the
surface decomposes as a union of cylinders.
A surface $S$ is {\em completely periodic} if any direction in which
there is at least one cylinder or a saddle connection connecting a
zero to itself is completely periodic.

It can be shown that if a translation surface $S$ has genus $g \geq 1$, then the number of zeroes of $\omega$, counting multiplicity, is $2g-2$.  Let $\alpha$ be a partition of $2g-2$, that is a collection of positive integers $\alpha_1, \hdots, \alpha_k$ such that $\alpha_1 + \hdots + \alpha_k = 2g-2$.  Thus, if $\omega$ has $k$ distinct zeroes, and $\alpha_i$ is the multiplicity of one of the $k$ zeroes, then $\omega$ determines a partition of $2g-2$.  Now, fix a genus $g$ and a partition $\alpha$ of $2g-2$.  We let $\cH(\alpha)$ denote the space of translation surfaces $(M, \omega)$ where $M$ is a Riemann surface of genus $g$ and $\omega$ is an abelian differential on $M$ whose zeroes correspond in the way described above to the partition $\alpha$. For example, the set of translation surfaces of genus $2$ is stratified as
$\cH(1,1) \cup \cH(2)$, where in $\cH(1,1)$ the 1-form has two distinct simple 
zeroes, and in $\cH(2)$ it has a single double zero. 

There is an action of the group $SL(2, \RR)$ on each stratum $\cH(\alpha)$.  If we think of a translation surface $S$ as a two-manifold with transition functions which are translations and which has a finite number of conical singularities, and if $g \in SL(2, \RR)$, then $gS$ is the two manifold which has charts which are given by $g \circ f_i$ where $f_i$ is a chart of $S$. It is easily checked that the transition functions of $gS$ are still translations and the number of singularities and their total angles are preserved. Alternatively, if $S=P_1 \cup \hdots \cup P_n$, where each $P_i$ is a polygon in the plane $\RR^2$, then $gS=gP_1 \cup \hdots \cup gP_n$.  A translation surface is said to be {\em Veech} if its stabilizer under this action of $SL(2, \RR)$ is a lattice.  It is a theorem of Veech that Veech surfaces are completely periodic, but the converse is not true. 

For a curve $\gamma \subset S$, let $p(\gamma)$ denote
$\int_\gamma(\omega)$, where $\omega$ is the $1$-form defining the
translation surface structure. We identify $\RR^2$ with $\CC$. 
Define a translation surface $S$ to be {\em quadratic} if there exists $d > 0$
(not a perfect square) such that $p(H_1(S, \ZZ)) \subset \QQ(\sqrt{d})
\times \QQ(\sqrt{d})$. Any Veech surface in genus $2$ (perhaps after
rescaling) is known to be
quadratic. 

Let $C$ be a cylinder of periodic trajectories in the direction $v$. 
We use the term {\em  width} for the length of the periodic trajectory, and {\em height}
for the distance across the cylinder in the orthogonal direction. If
$w$ is any other direction, we can talk about the {\em twist} of $C$ along
$w$.  Let $v^{\perp}$ be the direction orthogonal to $v$ given by the
right hand rule, and
let $w'$ be the vector in the direction $w$ whose projection onto
$v^{\perp}$ has length equal to the height of the cylinder. Now,
suppose that the angle $\theta$ between $v$ and $w$ satisfies $0\leq
\theta \leq \pi/2$.  Then we can define the twist $t$ in the direction
$w$
to be the length of the projection of $w'$ onto $v$ modulo the width
of the cylinder.  If $\pi/2 < \theta \leq \pi$, then we define the
twist to be the width of the cylinder minus the length of the
projection of $w'$ onto $-v$ modulo the width of the cylinder.

\begin{thm}
\label{theorem:explicit}
Let $S$ be a Veech surface in $\cH(2)$ which cannot be rescaled so that $p(H_1(S, \ZZ)) \subset \QQ \times \QQ$.   
Let
$v$ be a completely periodic direction for $S$. Suppose 
the cylinder decompositon in the direction $v$ has two
cylinders. Let $w_1, w_2, h_1, h_2,
t_1, t_2$
be the widths, heights and twists of these two cylinders.  
The twists are measured along some direction $w$ which we assume
is also completely periodic. After
rescaling the surface, we may assume that
these variables are algebraic integers in $\QQ(\sqrt d)$ where $d$ is some
square-free, positive integer. Then the following equations are satisfied:
\begin{align}
\label{eq:w1:barh1:w2:barh2}
w_1 \bar h_1 &=-w_2 \bar h_2 \\
\label{eq:twists}
\bar {w_1}t_1 + \bar{w_2}t_2 + \bar{w_1}w_2 &= w_1\bar{t_1} + w_2\bar{t_2}+w_1\bar{w_2}, \qquad {0 \le t_1 < w_1, 0 \le t_2 < w_2}
\end{align}
(where the bar denotes conjugation in $\QQ(\sqrt d)$.)

Conversely, let $S$ be a surface in $\cH(2)$.
Suppose there exists a direction $v$ in  which there is a
cylinder decomposition for which the heights, widths and twists (along
some direction $w$) belong
to $\QQ(\sqrt{d})$ and satisfy
\eqref{eq:w1:barh1:w2:barh2} and \eqref{eq:twists}. Then $S$ is Veech. 
\end{thm}

(If $S \in \cH(2)$ is Veech and has a one cylinder decomposition in some direction $v$, then $p(H_1(S, \ZZ)) \subset \QQ \times \QQ$.)

We note that a different classification of Veech surfaces in genus $2$, formulated in terms of Jacobians with real multiplication, was obtained by C. McMullen in \cite{ctm:billiards}.  

Let $\cO_d^+$ denote the positive algebraic integers in
$\QQ(\sqrt{d})$. 
\begin{lem}
\label{lemma:finetely:many:solutions}
For fixed integers $c_1$ and $c_2$, the equations
\eqref{eq:w1:barh1:w2:barh2}, \eqref{eq:twists} and
\begin{equation}
\label{eq:w1:h1:w2:h2:c1:c2}
w_1h_1 + w_2h_2= 2(c_1 + c_2 \sqrt{d})
\end{equation}
have finitely many solutions in $\cO_d^+$, 
up to the action by the group of units. Let us
denote this finite number by $H(c_1, c_2)$. 
\end{lem}

Let $\Omega$ be a set containing one element for each 
Veech surface in $\cH(2)$, where we
identify any two surfaces in the same $GL(2,\RR)$ orbit.  Let $P$
denote the set of positive solutions in algebraic integers in
$\QQ(\sqrt d)$ to \eqref{eq:w1:barh1:w2:barh2}
and \eqref{eq:twists}. 
Then in view of Theorem~\eqref{theorem:explicit} there
exists a map $f_1: P \to \Omega$. 
Now if $S$ is a Veech surface, then there is a surface $S'$  in the 
$GL(2,\RR)$ orbit of $S$ such that $p(H_1(S'))$ is a primitive
sublattice of $\cO_d
\times \cO_d$ with no common factor. The area of $S'$ is independent
of the choice of $S'$. Let $f_2: \Omega \to \cO_d^+$ denote the map
which sends (the equivalence class of) $S$ to the area of $S'$. 
\begin{cor}{\bf (Classification of Veech surfaces in $\cH(2)$).}
\label{cor:classification}
The maps $f_1: P \to \Omega$  and $f_2:\Omega \to
\cO_d^+$  are finite-to-one. In fact, 
the cardinality of $(f_1 \circ f_2)^{-1}(2(c_1 +c_2 \sqrt{d}))$ is
$H(c_1,c_2)$. 
\end{cor}

We now state some results for $\cH(1,1)$. 
Suppose $v$ is a completely periodic direction for $S \in
\cH(1,1)$. Without loss of generality, we may assume that in the
direction $v$, $S$ decomposes into three cylinders, and for $1 \le i
\le 3$, let $w_i$,
$h_i$ and $t_i$ denote the widths, heights and twists. After
renumbering, we may assume that $w_3 = w_1 + w_2$.  Denote $s_1 = h_1 +
h_3$, $s_2 = h_2 + h_3$, $\tau_1 = t_1 + t_3$, $\tau_2 = t_2 + t_3$. 
For $c_1 + c_2 \sqrt{d} \in \cO_d^+$, consider the following
equations: 
\begin{align}
\label{eq:w1:bars1:w2:bars2}
w_1 \bar s_1 &=-w_2 \bar s_2 \\
\label{eq:tau:twists}
\bar w_1 \tau_1 + \bar w_2 \tau_2 + \bar w_1 w_2 &= w_1 \bar \tau_1 +
w_2 \bar \tau_2 + w_1 \bar w_2 , \qquad {0 \le t_1 < w_1, 0 \le t_2 < w_2} \\
\label{eq:w1:s1:w2:s2:c1:c2}
w_1 s_1 + w_2 s_2 &= 2(c_1 + c_2 \sqrt{d})
\end{align}

\begin{thm}
\label{theorem:H11}
Let $\cM(c_1,c_2)$ denote the set of surfaces in $\cH(1,1)$ which are
completely periodic and for each completely periodic 
direction $v$, the lengths heights
and twists of the cylinder decomposition in that direction
satisfy \eqref{eq:w1:bars1:w2:bars2}, \eqref{eq:tau:twists}
and \eqref{eq:w1:s1:w2:s2:c1:c2}. Then $\cM(c_1,c_2)$ is a closed
subset of dimension $3$ which is invariant under the $SL(2,\RR)$
action. It is non-empty provided \eqref{eq:w1:bars1:w2:bars2}
\eqref{eq:tau:twists} and \eqref{eq:w1:s1:w2:s2:c1:c2} have a solution
in $\cO_d^+$, and any Veech 
surface in $\cH(1,1)$ is contained in such a subset. 
\end{thm}

We also obtain a classification of the completely periodic surfaces in
genus $2$. In fact
\begin{thm}
\label{theorem:cp:genus:2}
In $\cH(2)$ every completely periodic surface is Veech. In $\cH(1,1)$,
a surface is completely periodic if it has a cylinder
decomposition in some direction $v$, such that the heights widths and
twists satisfy
(\ref{eq:w1:bars1:w2:bars2}) (\ref{eq:tau:twists}) and
(\ref{eq:w1:s1:w2:s2:c1:c2}). Conversely, if $S \in \cH(1,1)$ is
completely periodic then after rescaling, either 
$p(H_1(S)) \subset \QQ \times \QQ$ or 
$p(H_1(S)) \subset \QQ(\sqrt{d}) \times \QQ(\sqrt{d})$ for some
square-free $d > 0$. If $S$ cannot be rescaled so that 
$p(H_1(S)) \subset \QQ \times \QQ$ and 
if $v$ is any completely 
periodic direction, then (\ref{eq:w1:bars1:w2:bars2})
(\ref{eq:tau:twists}) and (\ref{eq:w1:s1:w2:s2:c1:c2}) hold, 
where the $w_i, h_i, t_i$ denote the widths, heights and twists of the
cylinder decomposition along $v$ (and the twists are measured 
along another completely periodic direction $v'$). 
\end{thm}

\medskip
The proofs of the results stated here use the $J$ invariant for translation
surfaces defined  by Kenyon and Smillie in their paper ``Billiards on
rational angled triangles'' as follows.  First if $P$ is a polygon in
$\RR^2$ with vertices $v_1, \hdots, v_n$ in counterclockwise order
about the boundary of $P$, then $J(P)= v_1 \wedge v_2 + \hdots +
v_{n-1} \wedge v_n + v_n \wedge v_1$.  Then, if $S$ is a translation
surface with a cellular decomposition into planar polygons $P_1 \cup
\hdots \cup P_k$, $J(S)= \Sigma_{i=1}^{k} J(P_i)$.  Furthermore, there
are two linear projections $J_{xx}$ and $J_{yy}$ from $\RR^2 \wedge_{\QQ} \RR^2$
to $\RR \wedge_{\QQ} \RR$ given on basis elements as

$$J_{xx}\left(\binom{a}{b} \wedge \binom{c}{d}\right)=a \wedge c$$

$$J_{yy}\left(\binom{a}{b} \wedge \binom{c}{d}\right)=b \wedge d$$

\noindent and one linear projection $J_{xy}$ from $\RR^2 \wedge_{\QQ} \RR^2$
to $\RR \otimes_{\QQ} \RR$ defined on basis elements as 

$$J_{xy}\left(\binom{a}{b} \wedge \binom{c}{d}\right)=a \otimes d - c \otimes  b.$$

Let $v$ be some direction on a translation surface
$S$. If $v$ is either the horizontal or vertical direction, $J_{vv}$ has already
been defined. Suppose that $v$ is neither the horizontal nor vertical
direction, and let the vector $\binom{1}{q}$ represent the direction $v$.  If $g$ is any
matrix in $SL(2,\RR)$ such that $g\binom{1}{q}=\binom{1}{0}$, then we define
$J_{vv}(S)=J_{yy}(gS)$. It can be shown that this quantity is
independent of the choice of matrix $g$.  Also, given a pair of
linearly independent 
directions $v, w$ on $S$, let $v'$, $w'$ be vectors in the directions
$v$ and $w$ so that there exists an element $g \in SL(2,\RR)$ such that 
 $gv'=\binom{1}{0}$ and $gw'=\binom{0}{1}$ and define
 $J_{vw}(S)=J_{xy}(gS)$.  (In the case that $v=\binom{0}{1}$ and
 $w=\binom{1}{0}$, define $J_{vw}(S)=J_{yx}(S)$ where
 $J_{yx}(\binom{a}{b} \wedge \binom{c}{d})=b \otimes c - d
 \otimes a$.) It can be shown that a different choice of $v',w'$ will
 multiply each of the components of each of the tensors comprising
 $J_{vw}$ by the same number.

An alternative formula for the Kenyon-Smillie invariant is:
\begin{displaymath}
J = 2\sum_{i=1}^g p(a_i) \wedge p(b_i)
\end{displaymath} where $a_{i}, b_{i}$ is a symplectic homology basis
for the surface $S$. 

It is proven in Kenyon and Smillie that if the horizontal direction
$x$ on a translation surface $S$ is completely periodic, then $J_{yy}=0$.

We define a direction $v$ to be {\em
  homological} if there exists $\lambda \in H_1(S)$ and $r \in \RR$
such that $r v = p(\lambda)$. Note that $\lambda$ is not assumed to be
simple. 

We say that a translation surface $S$ has {\bf Property X} if for any
homological direction $v$ on $S$, $J_{vv}=0$. I have proven the
following theorems regarding property $X$.

\begin{thm}
\label{theorem:X:to:CP}
 If a genus $2$ translation surface has property $X$, then it is
  completely periodic. \end{thm}

\begin{thm} 
\label{theorem:Veech:to:X}
If a translation surface is completely periodic, then it has property
  $X$. In particular, any Veech surface has Property $X$.\end{thm}

\begin{thm} 
\label{theorem:H2:X:equals:Veech}
If a surface in $\cH(2)$ has property $X$, then it is Veech. \end{thm}

\begin{prop}
\label{prop:c1:equals:dc4}
 Let $S$ be a quadratic translation surface which cannot be rescaled so that $p(H_1(S, \ZZ)) \subset \QQ \times \QQ$.  Let $v$
and $w$ be any two linearly independent 
directions such that $J_{vv} = J_{ww} = 0$. 
We may write $J_{vw}(S)=c_1(1 \otimes
1 ) + c_2 (1 \otimes \sqrt{d}
) + c_3(\sqrt{d} \otimes 1 ) + c_4(\sqrt{d} \otimes \sqrt{d})$ where
$c_i \in \QQ$.
Then $S$ has property $X$ if and only if $c_2=c_3$ and
$c_1=dc_4$.
\end{prop}

A straightforward calculation shows the following:
\begin{lem}
\label{lemma:coords:to:X}
Let $S \subset \cH(2)$ be a quadratic translation surface which cannot be rescaled so that $p(H_1(S, \ZZ)) \subset \QQ \times \QQ$.   
Assume that $S$ has a two cylinder
decomposition in the horizontal direction.  Let $w_1, w_2, h_1, h_2,
t_1, t_2$ be the widths, heights and twists of these two cylinders. 
Then $J_{xx} = 0$ if and only if \eqref{eq:twists} holds. Furthermore
the equations $c_1 = d c_4$ and $c_2 = c_3$ defining property $X$ are
equivalent to \eqref{eq:w1:barh1:w2:barh2}. 
\end{lem}

Now Theorem~\ref{theorem:explicit} follows from
Theorem~\ref{theorem:Veech:to:X},
Theorem~\ref{theorem:H2:X:equals:Veech},
Proposition~\ref{prop:c1:equals:dc4},
and Lemma~\ref{lemma:coords:to:X}. 
The proof of Theorem~\ref{theorem:H2:X:equals:Veech} uses
Theorem~\ref{theorem:X:to:CP} and 
Lemma~\ref{lemma:finetely:many:solutions}.
The proof of Theorem~\ref{theorem:H11} also uses
Theorem~\ref{theorem:X:to:CP} and Lemma~\ref{lemma:finetely:many:solutions}.

We say that a translation surface $S$ is {\bf hyperperiodic} if any
homological direction on $S$ is completely periodic. 
Clearly hyperperiodicity implies both complete periodicity and property
$X$. The following theorem was obtained in joint work with A.~Eskin
and M.~Boshernitzan (see \cite{hyperperiodic}). 
\begin{thm}
\label{theorem:hyperperiodic}
Any quadratic Veech surface is hyperperiodic. Also, in genus $2$,
hyperperiodicity, complete periodicity and property X are equivalent. 
\end{thm}

We note that Theorem~\ref{theorem:hyperperiodic} gives a complete
description of the periodic directions in e.g. the regular pentagon. 
We also note that hyperperiodicty remains invariant under the
deformation where the absolute homology is fixed and the relative
homology changes.

\bigskip
Department of Mathematics, University of Chicago, Chicago IL~60637,
USA; \\
kcalta@math.uchicago.edu

\end{document}